\documentclass[11pt]{amsart}
\usepackage[foot]{amsaddr}
\usepackage{amssymb}
\usepackage[polish,english]{babel}
\usepackage{lmodern}
\usepackage[utf8]{inputenc}
\usepackage[T1]{fontenc}
\usepackage{amsmath,amscd}
\usepackage{graphicx}
\usepackage{underscore}
\usepackage[normalem]{ulem} 

\newcommand{\bR}{{\mathbb R}}

\usepackage[hidelinks]{hyperref}
\newcommand{\doi}[1]{\\ \href{https://doi.org/#1}{\texttt{https://doi.org/#1}}}

\newtheorem{theorem}{Theorem}
\newtheorem*{maintheorem*}{Main Theorem}
\newtheorem{proposition}[theorem]{Proposition}
\newtheorem{lemma}[theorem]{Lemma}
\newtheorem{remark}[theorem]{Remark}

\newtheorem{definition}[theorem]{Definition}

\title{An absorbing set for the Chialvo map}

\author[P.\ Pilarczyk]{Pawe\l{} Pilarczyk$^{\diamond}$}
\address{\rm $^{\diamond}$ Pawe\l{} Pilarczyk is with Faculty of Applied Physics and Mathematics \& Digital Technologies Center,
Gda\'{n}sk University of Technology,
ul.\ Narutowicza 11/12,
80-233 Gda\'{n}sk.
\newline
Orcid: \href{https://orcid.org/0000-0003-0597-697X}{0000-0003-0597-697X}
\newline
\url{https://www.pawelpilarczyk.com/}
}

\author[G.\ Graff]{Grzegorz Graff$^*$}
\address{\rm $^*$ Grzegorz Graff is with Faculty of Applied Physics and Mathematics \& BioTechMed Center,
Gda\'{n}sk University of Technology,
ul.\ Narutowicza 11/12,
80-233 Gda\'{n}sk.
\newline
Orcid: \href{https://orcid.org/0000-0001-5670-5729}{0000-0001-5670-5729}
\newline
[Corresponding author]
\newline
\url{https://pg.edu.pl/en/p/grzegorz-graff-8745}
\vskip 0pt
\copyright{} 2024. This manuscript version is made available under the CC-BY-NC-ND 4.0 license \url{https://creativecommons.org/licenses/by-nc-nd/4.0/}.
\vskip 0pt
This is the authors' accepted manuscript (AAM) of the paper
published in \textit{Communications in Nonlinear Science and Numerical Simulation},
Volume 132, May 2024, 107947 \doi{10.1016/j.cnsns.2024.107947}
}

\date{}

\begin{document}

\begin{abstract}
The classical Chialvo model, introduced in 1995, is one of the most important models that describe single neuron dynamics. In order to conduct effective numerical analysis of this model, it is necessary to obtain a rigorous estimate for the maximal bounded invariant set. We discuss this problem, and we correct and improve the results obtained by Courbage and Nekorkin [Internat. J.\ Bifur.\ Chaos Appl.\ Sci.\ Engrg. 20 (2010), 1631--1651.] In particular, we provide an explicit formula for an absorbing set for the Chialvo neuron model. We also introduce the notion of a weakly absorbing set, outline the methodology for its construction, and show its advantage over an absorbing set by means of numerical computations.

\smallskip

\noindent \textbf{\keywordsname:} Chialvo map, neuron model, grid, rigorous numerics, absorbing set.
\end{abstract}


\maketitle


\section{Introduction}
\label{sec:intro}

Studying neural dynamics is nowadays one of the most important topics in mathematical biology. In spite of considerable amount of research that has been conducted in this area for neural networks so far, understanding neural processes at the level of individual cells is still an important scientific challenge.

In this paper, we focus on one of the simplest yet mathematically demanding models:
the classical two-dimensional discrete-time model of a single neuron that was proposed by Dante R.\ Chialvo in 1995 \cite{Chialvo}:
\begin{equation}
\label{eq:main}
\left\{\begin{array}{rcl}
x_{n+1} & = & x_n^2 \exp(y_n-x_n) + k,
\\
y_{n+1} & = & ay_n -bx_n + c.
\end{array}
\right.
\end{equation}

In this model, the variable $x$ (also called \emph{activation}) is related to the evolution of the membrane potential of a neuron, and $y$ is a recovery-like variable. 
The Chialvo model represents so-called excitable dynamics, with the variable $y$ responsible for fast recovery mechanisms.

Let us mention the fact that the equations (\ref{eq:main}) introduced by Chialvo in \cite{Chialvo} constitute one of the first two-dimensional discrete-time models of a single neuron, which was later followed by many similar models, such as the Rulkov model, the Courbage–Nekorkin–Vdovin model, or the discrete Izhikevich model; see \cite{izi2003} and the review papers \cite{courbage2010,Girardi,ibarz2011}.

Although the Chialvo equations yield a considerably simplified model of the real neuron, they capture all the important key features of neuronal dynamics. In particular, the presence of a complicated structure of bifurcation patterns, multiple attractors of various kinds, and chaotic behavior of the system observed within wide ranges of parameters mirror the rich complexity of the actual behavior of real neurons.

Since its introduction, the Chialvo model has been subject of growing interest. Many authors were studying various aspects of the model, such as the presence of the so-called comb-shaped chaotic structure~\cite{NewPaperOnChialvo}, different types of bifurcations~\cite{Jing,chialvo2d}, the role of switching mechanism for the dynamics of the model~\cite{Yang}, reduction to one-dimensional
version \cite{Trujillo}, stochastic effects of neural activity~\cite{Bashkirtseva}, or dynamical effects of electromagnetic flux on the model~\cite{chialvoMagneto}. The Chialvo map is also used to define local dynamics of nodes that form a network used to model synchronization of many neurons~\cite{Kamal}. Thorough understanding of the Chialvo model itself is thus important for understanding the dynamics of neural networks.

The model has four real parameters: $a>0$ -- the time constant of the recovery, $b>0$ -- the activation dependence of the recovery process, $c$ -- the offset, and $k$ -- an additive perturbation.

A systematic analysis of the model by the use of rigorous numerical methods has been conducted in \cite{chialvo2d} for a wide parameter range that covers the classical Chialvo considerations~\cite{Chialvo}. Namely, the values of $a = 0.89$ and $c = 0.28$ were fixed, and the two other parameters were varied in the range $(b,k) \in  [0,1] \times [0,0.2]$.
This range was partitioned into classes corresponding to qualitatively different types of dynamics (continuation classes). Moreover, bifurcation patterns across the spectrum of considered parameters were described, and insight into the dynamics within each continuation class was given. This information was conveyed through the utilization of the Conley index and Morse decomposition. Additional analysis was conducted aimed at determining how alterations in parameter values were inducing changes in the system's dynamics, with emphasis on detection of chaotic dynamics.

An important question when studying dynamics of such systems numerically is to identify a bounded set $P$ such that all the recurrent dynamics of the system is contained in $P$ (this set may depend on the actual values of parameters of the system considered).
It would be ideal if $P$ absorbs every forward trajectory, as we define below.
\begin{definition}[absorbing set]
\label{def:AbsorbingSet}
A set $P \subset X$ is called an \emph{absorbing set} for a map $f \colon X \to X$ if for every $x \in X$ there exists $N_0 \geq 0$ such that $f^n (x) \in P$ for all $n \geq N_0$.
\end{definition}
\noindent
In this case, no trajectory would leave $P$ in forward time, and numerical computations can be restricted to the set $P$ without losing any important long-term dynamical features.

However, a full picture of all the recurrent dynamics can also be reconstructed if we find a set $P$ that is merely crossed by every forward trajectory in the system. This kind of a set may be defined as follows.
\begin{definition}[weakly absorbing set]
\label{def:WeaklyAbsorbingSet}
A set $P \subset X$ is called a \emph{weakly absorbing set} for a map $f \colon X \to X$ if for every $x \in X$ there exists some $n \geq 0$ such that $f^n (x) \in P$.
\end{definition}
\noindent
Note that given a weakly absorbing set $P$, the union of its forward iterates $\bigcup_{n \geq 0} f^n(P)$ is an absorbing set. With weaker properties required, it may be expected that such a set will be smaller and thus easier to tackle numerically; see Section~\ref{sec:D1plusHuge} for some comparison.

For the purpose of the potential use of methods of interval arithmetic in studying the dynamics (such as, for example, the methods used in \cite{Arai,GAIO,Knipl,Miyaji,chialvo2d}), $P$ should be of the form of the Cartesian product of two bounded intervals $P = I_x \times I_y$.
Then one could subdivide such a set $P$ into a finite collection of compact rectangles with pairwise disjoint interiors by introducing a grid defined by subdividing $I_x$ and $I_y$ into intervals: $x_0 < x_1 < \cdots < x_N$ and $y_0 < y_1 < \cdots <y_M$, where $I_x = [x_0, x_N]$ and $I_y = [y_0, y_M]$. Each grid element is of the form $[x_i,x_{i+1}] \times [y_i,y_{j+1}]$. Computing images of such grid elements directly using interval arithmetic makes it possible to conduct rigorous analysis of dynamics at the finite resolution, and, for example, to use fast algorithms to split $P$ into recurrent subsets of various kinds. In particular, one can use algorithms introduced in \cite{MMP2005,PS2008} to compute Morse decompositions and the Conley index of isolating neighborhoods built of grid elements. The assumption that $P$ contains all the recurrent dynamics of the system guarantees in this context that no recurrent dynamics has been missed and allows treating the results of these computations as computer assisted proof.

In this paper, we study the problem of determining an absorbing set and a weakly absorbing set for the Chialvo map~\eqref{eq:main}, taking as a starting point the construction from \cite{courbage2010} that we correct and modify.
Our main results (Theorems \ref{thm:absorbing} and~\ref{thm:weak}) can be summarized as follows.

\begin{maintheorem*}
Let $a \in (0,1), b \in (0,1), c \in (0,1), k>0$. Consider
\begin{eqnarray*}
\hat{v} & := & (c-bk)/(1-a), \qquad \text{see \eqref{eq:v}} \\
\hat{u}_1 & := & 4 / e^2 \cdot \exp(\hat{v}) + k, \qquad \text{see \eqref{eq:u1}} \\
\hat{w}_1 & := & (c - b \hat{u}_1) / (1-a), \qquad \text{see \eqref{eq:w1}}
\end{eqnarray*}
Then the set $\hat{D}_1^+ := [k, \hat{u}_1] \times [\hat{w}_1, \hat{v}]$ is an absorbing set for the Chialvo map~\eqref{eq:main}.

Moreover, if $a+b \geq 1$ and $\hat{u}_2$ is defined as described in Section~\ref{sec:C2}, see~\eqref{eq:u2}, and
\begin{eqnarray*}
\hat{w}_2 & := & (c - b \hat{u}_2) / (1-a), \qquad \text{see \eqref{eq:w2}}
\end{eqnarray*}
then the set $\hat{D}_2^+ := [k, \hat{u}_2] \times [\hat{w}_2, \hat{v}]$ is a weakly absorbing set for the Chialvo map~\eqref{eq:main}.
\end{maintheorem*}

The paper is organized in the following way. In Section \ref{sec:problem}, we state the problem and we describe the result of \cite{courbage2010} that unfortunately turned out to be incorrect. In the next section, we show counterexamples to some statements that appear in the construction of an absorbing set proposed in \cite{courbage2010}. Our main results are described in the next two sections: In Section \ref{sec:abs}, we provide the construction of an absorbing set for the Chialvo map~\eqref{eq:main}, which is summarized in Theorem~\ref{thm:absorbing}. Section~\ref{sec:weak} is devoted to the construction of a weakly absorbing set, with the main result stated in Theorem~\ref{thm:weak}. Finally, in Section~\ref{sec:num}, we conduct numerical study of various features of the absorbing and weakly absorbing sets; in particular, we compare their size, and we show how they approximate the global attractor observed numerically for some classical parameter values.
 

\section{Statement of the problem}
\label{sec:problem}

When studying numerically the dynamics of the Chialvo model, one needs to have a bounded set $P \subset \bR^2$ that captures all the recurrent dynamics. Ideally, $P$ would be an absorbing set for the Chialvo map, as in Definition~\ref{def:AbsorbingSet}.

Let us remark at this point that sometimes one might be interested in a local absorbing set. This is a set that absorbs trajectories from its certain neighborhood (cf., for example,~\cite{Maistrenko1}).
For the reasons described earlier, however, we will seek a globally absorbing set $P$.
As explained above, we would like to have $P$ in the form of the Cartesian product of intervals: $P = I_x \times I_y$. Obviously, the intervals that define $P$ would depend on the parameters of the system, and we are interested in a relatively tight estimate.
In what follows, we consider the following wide ranges of the parameters of the Chialvo map
\begin{equation}
\label{eq:ranges}
a \in (0,1), \quad b \in (0,1), \quad c \in (0,1), \quad k>0.
\end{equation}

In an interesting survey paper \cite{courbage2010}, the authors made an attempt to identify an absorbing set for the Chialvo model. Namely, on page 1643, after (30) of~\cite{courbage2010}, the following set $D^+$ is defined that is claimed to be an absorbing set:
\begin{equation}
    \label{AbsorbingSet}
    D^+:=\left\{(x,y): \; k<x<\frac{2a}{a-b}, \; \frac{c(a-b)-2ab}{1-a}<y<\frac{c-bk}{1-a}\right\}.
\end{equation}

Although the general idea of finding such a set used in \cite{courbage2010} looks promising, the reasoning in Section 3.2.1 of that paper seems to have some gaps and omissions, and the set $D^+$ turns out not to be absorbing.

\begin{figure}[htbp]
\includegraphics[width=1.0\textwidth]{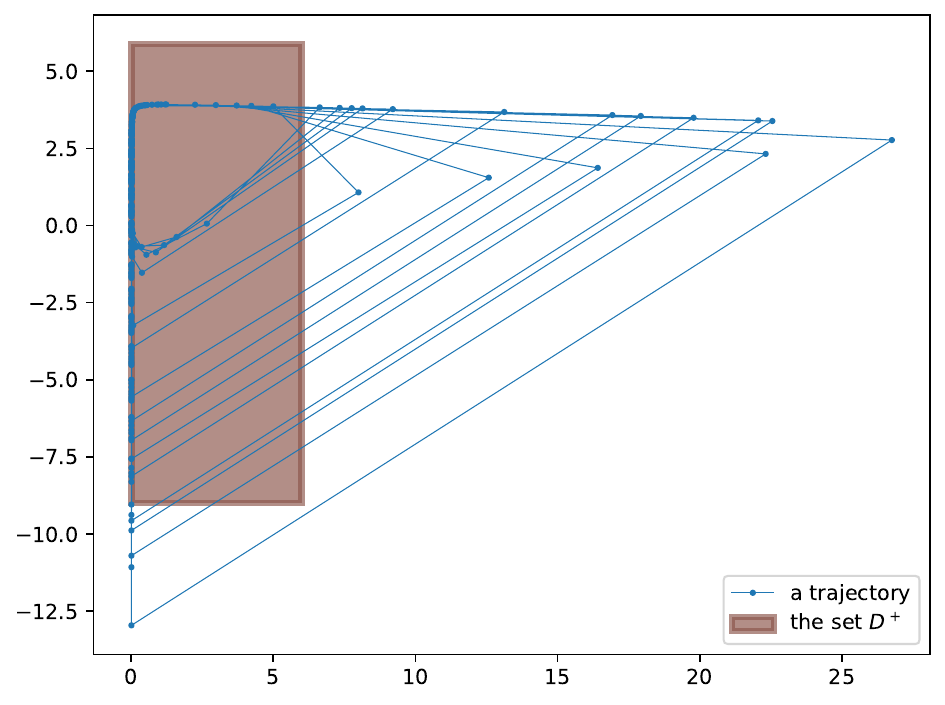}
\caption{\label{fig:badtrajectory}%
Three hundred iterates of a seemingly chaotic trajectory observed in the Chialvo model for $a=0.9$, $b=0.6$, $c=0.6$, and $k=0.02$ after initial $100{,}000$ iterates that started at $(x_0,y_0)=(1,1)$. Consecutive points on the trajectory are joined by thin straight lines. The set $D^+$ calculated in \eqref{eq:absorbing1} according to \eqref{AbsorbingSet} is shown in brown.}
\end{figure}

For example, consider the Chialvo model with $a=0.9$, $b=0.6$, $c=0.6$, and $k=0.02$.
Formula \eqref{AbsorbingSet}
yields the following bounds on the absorbing set:
\begin{equation}
\label{eq:absorbing1}
x \in [0.02,6], y \in [-9,5.88].
\end{equation}
However, numerical simulations show that a trajectory that starts from the initial condition $(x,y)=(1,1)$, after $100{,}000$ iterations (when it should have already settled down on an attractor) is not contained in $D^+$ defined by \eqref{eq:absorbing1}. Its $300$ next iterates are shown in Figure~\ref{fig:badtrajectory} along with the set $D^+$.

Although the trajectory shown in Figure~\ref{fig:badtrajectory} is not absorbed by $D^+$, it nevertheless crosses $D^+$ several times. Therefore, one might be tempted to suspect that maybe $D^+$ is weakly absorbing in the sense of Definition~\ref{def:WeaklyAbsorbingSet}.
Unfortunately, as we show in Section~\ref{sec:details}, the set $D^+$ defined by \eqref{AbsorbingSet} is not even weakly absorbing.

The aim of the paper is to provide counterexamples to the statements in \cite[Section 3.2.1]{courbage2010} as well as to introduce a new reasoning and formulas that allow one to define a correct absorbing set for the Chialvo map.


\section{Counterexamples to the statements about $D^+$}
\label{sec:details}

\begin{figure}[htbp]
\includegraphics[width=0.5\textwidth]{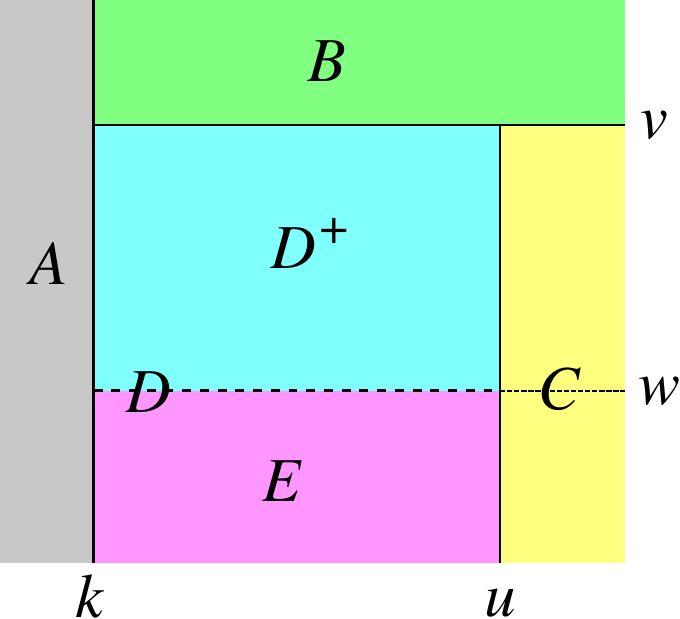}
\caption{\label{fig:regions}%
Regions in the phase space discussed in the paper. The actual values of $u$, $v$ and $w$ are indicated in the text.}
\end{figure}

Let us discuss the analysis of the Chialvo model conducted in \cite[Section 3.2.1]{courbage2010} regarding the absorbing region $D^+$ that is supposed to attract all trajectories with initial conditions outside it. See Figure~\ref{fig:regions} for reference to the regions discussed.

The regions discussed in \cite{courbage2010} are defined (up to the boundary) as follows:
\begin{eqnarray}
\label{eq:A} A & = & \{ x < k \}, \\
\label{eq:B} B & = & \{ x \geq k,\; y > v \}, \\
\label{eq:C} C & = & \{ x > u,\; y \leq v \}, \\
\label{eq:D} D & = & \{ k \leq x \leq u,\; y \leq v \}, \\
\label{eq:E} E & = & \{ k \leq x \leq u,\; y < w \} \subset D, \\
\label{eq:Dplus}  D^+ & = & D \setminus E,
\end{eqnarray}
where the boundaries between the regions are at the following levels:
\begin{eqnarray}
\label{eq:u0} u & = & \frac{2a}{a-b}, \\
\label{eq:v0} v & = & \frac{c-bk}{1-a}, \\
\label{eq:w0} w & = & \frac{c(a-b)-2ab}{1-a}.
\end{eqnarray}

In order to prove that $D^+$ is an absorbing region, the following reasoning was conducted in \cite{courbage2010}.

Assume that $a,b,c \in (0,1)$ and $k>0$.
Denote $(\bar{x},\bar{y}) := f(x,y)$.

If $(x,y) \in A$ then it follows from the first equation of \eqref{eq:main} that $\bar{x} \geq k$, and thus any trajectory leaves the region $A$ in just one iterate. From now on it is assumed that $x \geq k$.

If $y > v$ then from the second equation  of \eqref{eq:main} and from the fact that $x \geq k$ it follows (after a simple calculation) that $\bar{y} < y$. Note that obviously $\bar{x} > k$. As a consequence, any trajectory that starts in $B$ tends toward its lower boundary, eventually leaves $B$ after a finite number of iterates (this is not immediate, but we later provide proof of this fact; see Lemma~\ref{lem:downB} and Proposition~\ref{prop:B}), and enters the region $\{ x > k,\; y < v \}$.

The authors of \cite{courbage2010} prove that $f$ is one-to-one in the region $D$, which is a valuable observation, but does not seem relevant to the fact that $D^+$ is to be an absorbing region.

Then in \cite{courbage2010} it is claimed that $\bar{x} < x$ for every $(x,y) \in D$. Let us call this ``Statement~1.''

Right afterwards there is a claim that every trajectory that starts in Region $C$ (see Figure~\ref{fig:regions}) tends to Region $D$, but this is not proven. Let us call this ``Statement~2.''

Furthermore, it is claimed that if $(x,y) \in E$ (see Figure~\ref{fig:regions}) then $\bar{y} > y$. Let us call this ``Statement~3.'' Together with Statement~1 (saying that $\bar{x} < x$ for every $(x,y) \in D \supset E$) this is claimed in \cite{courbage2010} to prove that any trajectory starting in $E$ moves upward (in the direction of $D^+$) and does not leave the belt $k < x < u$ as long as it is below $w$.

Then it is claimed in \cite{courbage2010} that, as a consequence of the above, Region $D^+$ ``attracts all trajectories with initial conditions outside it,'' and thus $D^+$ is attracting as in Definition~\ref{def:AbsorbingSet}.

Unfortunately, the reasoning in \cite{courbage2010} has some gaps. First, the three statements (Statement~1, Statement~2 and Statement~3) are not true. 

Second, it turns out that the set $D^+$ is not absorbing  and not even weakly absorbing.

Let us now show specific counterexamples for the three statements discussed above.

Counterexample for Statement~1:
\begin{equation}
    x=0.2,\; y=-0.4,\; a=0.2,\; b=0.1,\; c=0.1,\; k=0.179.
\end{equation}
In this case $\bar{x}=0.201$ is not smaller than $x=0.2$. Note that $u=4$ and $v=0.103$.

Counterexample for Statement~2:
\begin{equation}
    x=2.6,\; y=1.6,\; a=0.5,\; b=0.1,\; c=0.9,\; k=0.116.
\end{equation}
In this case $\bar{x}=2.602$ is not smaller than $x=2.6$. Note that $u=2.5$ and $v=1.777$.

Counterexample for Statement~3:
\begin{equation}
    x=2,\; y=-0.4,\; a=0.3,\; b=0.2,\; c=0.1,\; k=0.01.
\end{equation}
In this case $\bar{y}=-0.42$ is not greater than $y=-0.4$. Note that $u=6$ and $w=-0.157$.

In fact, it turns out that the set $D^+$ is not absorbing (see Definition~\ref{def:AbsorbingSet}) and not even weakly absorbing (see Definition~\ref{def:WeaklyAbsorbingSet}). Namely, some trajectories never enter it in forward time. For example, one can find an attracting fixed point outside Region~$D^+$ for some parameter values. Here is one:
\begin{equation}
\label{eq:counter3}
    a=0.6,\; b=0.1,\; c=0.9,\; k=0.01,\; (x,y)=(2.553,1.612).
\end{equation}
Note that in this case $u=2.4$, $v=2.248$, and $w=0.825$. This fixed point is precisely in the interior of Region~$C$ whose analysis was omitted in~\cite{courbage2010}.

Additionally, we would like to point out the fact that if $a < b$ then we have $u < 0$ and thus $D^+ = \emptyset$ when $k > 0$. Moreover, if $a = b$ then the upper bound on $x$ given by $v = \frac{2a}{a-b}$ yields $+\infty$ and is thus useless.


\section{New construction of an absorbing set for the Chialvo map}
\label{sec:abs}

In order to define a correct absorbing set for the Chialvo map, we use the same regions as shown in Figure~\ref{fig:regions}, except we redefine the values of $u$, $v$ and $w$, which we now denote $\hat{u}$, $\hat{v}$ and $\hat{w}$, respectively.
The counterparts of the regions $A$, $B$, $C$, $D$, $D^+$ and $E$ will be denoted as $\hat{A}$, $\hat{B}$, $\hat{C}$, $\hat{D}$, $\hat{D}^+$ and $\hat{E}$, respectively.
Like previously, let us denote $(\bar{x},\bar{y}) := f(x,y)$, and additionally let us denote $(\bar{\bar{x}},\bar{\bar{y}}) := f(\bar{x},\bar{y})$.

The purpose of the definitions and statements in this section is to ultimately define the set $\hat{D}^+$ in such a way that it is absorbing, or at least weakly absorbing. In order to achieve this, we define the regions $\hat{A}$, $\hat{B}$, $\hat{C}$ and $\hat{E}$ around $\hat{D}^+$ and prove that every trajectory leaves these regions in forward time. If the trajectories leave these regions permanently then they enter $\hat{D}^+$ and stay there forever; then the set $\hat{D}^+$ is absorbing; see Theorem~\ref{thm:absorbing}. Otherwise, if every trajectory merely crosses $\hat{D}^+$ in forward time but may leave $\hat{D}^+$ later then the set is weakly absorbing; see Theorem~\ref{thm:weak}.


\subsection{\textbf{Region $\hat{A}$.}}
We begin with the same first step as in \cite{courbage2010}, and we define
\begin{equation}
\label{eq:regionA}
\hat{A} := \{ x < k \} = A.
\end{equation}
Notice that for every $(x,y) \in \bR^2$, we have $\bar{x} \geq k$, and if $x \geq k$ then $\bar{x} > k$; in particular, the image of every point in Region~$\hat{A}$ is situated in the closed half-plane bounded by the line $\{x=k\}$ from the left. We summarize these findings in the following statement.
\begin{proposition}[Region $\hat{A}$]
\label{prop:A}
For every $(x,y) \in \bR^2$, we have $(\bar{x},\bar{y}) \notin \hat{A}$. Moreover, if $(x,y) \notin \hat{A}$ then $\bar{x} > k$.
\end{proposition}
As a consequence, from now on we are only interested in trajectories that start to the right of the line $\{x = k\}$.


\subsection{\textbf{Region $\hat{B}$.}}
\label{sec:regionB}
Let us define
\begin{equation}
\label{eq:v}
\hat{v} := \frac{c - bk}{1 - a}
\end{equation}
as in \cite{courbage2010}, and then let us set
\begin{equation}
\label{eq:regionB}
\hat{B} := \{ k \leq x,\; \hat{v} < y \} = B.
\end{equation}

Whenever $(x,y) \in \hat{B}$, that is, $x > k$ and $y > \hat{v}$, we have $\bar{y} < y$ (see~\cite{courbage2010}), which can be checked with the following calculation:
\begin{equation}
\label{eq:yyBweak}
\bar{y} - y = -\underbrace{(1-a)}_{>0} y - b k + c < -(1-a)\frac{c-bk}{1-a} + (c-bk) = 0.
\end{equation}
This implies the fact that every trajectory moves downward in Region~$\hat{B}$. The following lemma provides a lower bound on the ``speed'' of this movement, which will be used in Proposition~\ref{prop:B} to prove that it will actually leave $\hat{B}$.
\begin{lemma}[moving downward in Region $\hat{B}$]
\label{lem:downB}
There exists a positive number $s > 0$ such that for every $(x,y) \in \hat{B}$ at least one of the following three conditions holds: $y - \bar{y} \geq s$, $(\bar{x},\bar{y}) \notin \hat{B}$, or $\bar{y} - \bar{\bar{y}} \geq s$.
\end{lemma}
\begin{proof}
Let us take any $(x,y) \in \hat{B}$ and use the fact that $y > \hat{v}$:
\begin{multline}
\label{eq:yyB}
y - \bar{y} = y - ay + bx - c =
(1-a)y + bx - bk + bk - c \\
\geq (1-a) \frac{c-bk}{1-a} - (c-bk) + b(x-k) =
b(x-k).
\end{multline}
Take $\varepsilon := \min \{k^2 e^{\hat{v}-k-1}, 1\}$.
If $x \geq k + \varepsilon$ then it follows from \eqref{eq:yyB} that $y - \bar{y} \geq b \varepsilon$.
Otherwise, $x < k + \varepsilon \leq k + 1$. Recall that $x \geq k$ and $y > \hat{v}$, and calculate:
\begin{equation}
\bar{x} = x^2 e^{y-x} + k > k^2 e^{\hat{v}-k-1} + k \geq \varepsilon + k,
\end{equation}
and one can use the argument above to show that $\bar{y} - \bar{\bar{y}} \geq b\varepsilon$, unless $(\bar{x},\bar{y}) \notin \hat{B}$.
As a consequence, one can take $s := b \varepsilon$, which completes the proof.
\end{proof}

The next lemma shows that once a trajectory has left Region~$\hat{B}$, it cannot come back.
\begin{lemma}[no return to Region $\hat{B}$]
\label{lem:noB}
If $(x,y) \notin \hat{A} \cup \hat{B}$ then $(\bar{x},\bar{y}) \notin \hat{B}$.
\end{lemma}
\begin{proof}
If $x>k$ and $y<\hat{v}$ then we have the following:
\begin{multline}
\bar{y}= a y - b x + c < a y - bk + c \\
< a\frac{c-bk}{1-a}-bk+c = a\frac{c-bk}{1-a} + \frac{1-a}{1-a} (c-bk) = \frac{c-bk}{1-a}=\hat{v}.
\end{multline}
This shows that $\bar{y} < \hat{v}$, which completes the proof.
\end{proof}

The key properties of Region~$\hat{B}$ are summarized in the following result that is a consequence of \eqref{eq:yyBweak} combined with Lemmas \ref{lem:downB} and~\ref{lem:noB}.
\begin{proposition}[Region $\hat{B}$]
\label{prop:B}
Every trajectory starting in $\hat{B}$ leaves the set $\hat{B}$ after a finite number of iterations. Moreover, if $(x,y) \notin \hat{A} \cup \hat{B}$ then $(\bar{x},\bar{y}) \notin \hat{B}$.
\end{proposition}


\subsection{\textbf{Regions $\hat{C}_1$ and $\hat{D}_1$.}}
\label{sec:C1}

We are going to define two counterparts $\hat{C}_1$ and $\hat{C}_2$ of Region $C$ defined in Section~\ref{sec:details} in the form $\hat{C}_i = \{ x > \hat{u}_i,\; y \leq \hat{v} \}$, $i=1,2$. Region~$\hat{C}_1$ is defined below, and we define Region~$\hat{C}_2$ in Section~\ref{sec:C2}. Since Regions $D$, $D^+$ and $E$ depend on $\hat{u}$, we will also have two counterparts of each of these regions.

Region $\hat{C}_1$ will be defined by the appropriate choice of the value of $\hat{u}_1$ in such a way that after just one iteration the points from $\hat{C}_1$ are transferred to the area
\begin{equation}
\label{eq:regionD1}
\hat{D}_1 := \{k \leq x \leq \hat{u}_1,\; y \leq \hat{v}\}.
\end{equation}

Recall that
\begin{equation}\label{mainx}
\bar{x} = x^2 e^{y-x} + k.
\end{equation}
Consider the function
\[
g(x)=\frac{x^2}{e^x} e^{\hat{v}} + k.
\]
By computing the derivative and analyzing its sign, one can see that $g$ has a global maximum at $x=2$:
\begin{equation}
\label{eq:gmax}
g(x) \leq g(2) = \frac{4}{e^2}e^{\hat{v}} + k.
\end{equation}
Let us define
\begin{equation}
\label{eq:u1}
\hat{u}_1 := \frac{4}{e^2}e^{\hat{v}} + k =
\frac{4}{e^2}e^{(c - bk) / (1 - a)} + k,
\end{equation}
and thus we define Region~$\hat{C}_1$ as follows:
\begin{equation}
\label{eq:regionC1}
\hat{C}_1 := \{\hat{u}_1 < x,\; y \leq \hat{v}\}.
\end{equation}
The above reasoning implies the following key property of Regions $\hat{C}_1$ and~$\hat{D}_1$:
\begin{proposition}[Regions $\hat{C}_1$ and $\hat{D}_1$]
\label{prop:C1}
Let $(x,y) \in \hat{C}_1 \cup \hat{D}_1$. Then $(\bar{x},\bar{y}) \in \hat{D}_1$.
\end{proposition}


\subsection{\textbf{Region $\hat{E}_1$.}}
\label{sec:E1}
Consider Region~$\hat{E}_1$ of the form
\begin{equation}
\label{eq:regionE1}
\hat{E}_1 =\{ k \leq x \leq \hat{u}_1,\; y < \hat{w}_1 \} \subset \hat{D}_1,
\end{equation}
where $\hat{w}_1 \leq \hat{v}$.
We would like to determine a possibly high value of $\hat{w}_1$ for which $y < \hat{w}_1$ implies $\bar{y} > y$, that is, all trajectories in Region~$\hat{E}_1$ move upward. For that purpose, we ask how small $y$ has to be so that $\bar{y} - y > 0$.
We propose the following:
\begin{equation}
\label{eq:w1}
\hat{w}_1 := \frac{c - b \hat{u}_1}{1-a}.
\end{equation}
To see that this is a good choice, take any point $(x,y) \in \hat{E}_1$. Then $y < \hat{w}_1$ and $x \leq \hat{u}_1$, and we can calculate as follows:
\begin{multline}
\label{eq:yyE}
\bar{y} - y = -(1-a) y - b x + c \\
> -(1-a) \hat{w}_1 - b \hat{u}_1 + c =
-(1-a) \frac{c - b \hat{u}_1}{1-a} + c - b \hat{u}_1 = 0.
\end{multline}

Similarly as it was the case of Region~$\hat{B}$ considered in Section~\ref{sec:regionB}, the inequality \eqref{eq:yyE} is not enough
to prove that a trajectory that starts in Region~$\hat{E}_1$ actually leaves this region in a finite number of steps. Therefore, similarly to Lemma~\ref{lem:downB}, we need the following result:

\begin{lemma}[moving upward in Region $\hat{E}_1$]
\label{lem:upE1}
There exists a positive number $s > 0$ such that for every $(x,y) \in \hat{E}_1$ at least one of the following three conditions holds: $\bar{y} - y \geq s$, $(\bar{x},\bar{y}) \notin \hat{E}_1$, or $\bar{\bar{y}} - \bar{y} \geq s$.
\end{lemma}
\begin{proof}
Take any point $(x,y) \in \hat{E}_1$. Use the fact that $y < \hat{w}_1$ and calculate:
\begin{multline}
\bar{y} - y = -(1-a) y - bx + c \\
\geq -(1-a) \frac{c - b\hat{u}_1}{1-a} - b x + c = b (\hat{u}_1 - x).
\end{multline}

Take any $\varepsilon > 0$ such that $\hat{u}_1 - \varepsilon > k$ and $\varepsilon \leq \tfrac{4}{e^2}(e^{\hat{v}} - e^{\hat{w}_1})$.
If $x < \hat{u}_1 - \varepsilon$ then $\bar{y} - y \geq b \varepsilon$.
Otherwise, $x \in [\hat{u}_1 - \varepsilon, \hat{u}_1]$. We shall show that $\bar{x} < \hat{u}_1 - \varepsilon$, so the inequality will hold at the next iteration, unless $(\bar{x},\bar{y}) \notin \bar{E}_1$. We use the fact that $y \leq \hat{w}_1$ and the definition \eqref{eq:u1} of $\hat{u}_1$:
\begin{multline}
\hat{u}_1 - \bar{x} = \frac{4}{e^2} e^{\hat{v}} + k - x^2 e^{y-x} - k \\
\geq \frac{4}{e^2} e^{\hat{v}} - \underbrace{\frac{x^2}{e^x}}_{\leq 4 / e^2} e^{\hat{w}_1} \geq  \frac{4}{e^2} (e^{\hat{v}} - e^{\hat{w}_1}) \geq \varepsilon.
\end{multline}
One can therefore take $s := b \varepsilon$, which ends the proof.
\end{proof}

The next lemma shows that once a trajectory has left Region~$\hat{E}_1$, it cannot come back.
\begin{lemma}[no return to Region $\hat{E}_1$]
\label{lem:noE1}
If $(x,y) \in \hat{D}_1 \setminus \hat{E}_1$ then $(\bar{x},\bar{y}) \notin \hat{E}_1$.
\end{lemma}
\begin{proof}
Take any $(x,y) \in \hat{D}_1 \setminus \hat{E}_1$.
Under this assumption, $y \geq \hat{w}_1$. Let us calculate a lower bound on $\bar{y}$:
\begin{multline}
\bar{y} = ay - bx + c \geq a \hat{w}_1 - b \hat{u}_1 + c \\
= a \frac{c - b \hat{u}_1}{1 - a} + \frac{1-a}{1-a} (c - b \hat{u}_1) = \frac{c - b \hat{u}_1}{1 - a} = \hat{w}_1.
\end{multline}
Therefore, $(\bar{x},\bar{y}) \notin \hat{E}_1$ indeed.
\end{proof}

To sum up, we have the following features of Region~$\hat{E}_1$:
\begin{proposition}[Region $\hat{E}_1$]
\label{prop:E1}
Every trajectory that starts in the set $\hat{E}_1$
leaves $\hat{E}_1$ in a finite number of iterations.
Moreover, if $(x,y) \in \hat{D}_1 \setminus \hat{E}_1$
then $(\bar{x},\bar{y}) \notin \hat{E}_1$.
\end{proposition}

\begin{remark}
For the parameters in counterexample \eqref{eq:counter3}, we would have $\hat{w}_1 = 1.374$, which is below the $y$ coordinate of the fixed point $(x,y)=(2.553,1.612)$ found for those parameter values in numerical simulations. Therefore, the fixed point is not in Region~$\hat{E}_1$.
\end{remark}


\subsection{\textbf{Region $\hat{D}_1^+$.}}
\label{sec:D1plus}

Let us now gather the properties of Region $\hat{D}_1^+$ bounded by Regions $\hat{A}$, $\hat{B}$, $\hat{C_1}$ and $\hat{E_1}$:
\begin{equation}
\label{eq:d1plus}
\hat{D}_1^+ := [k, \hat{u}_1] \times [\hat{w}_1,\hat{v}]
\end{equation}

We first notice that the obtained set is non-empty for all the ranges of parameters $a$, $b$, $c$ and $k$ that we consider.

\begin{proposition}[non-emptiness of Region $\hat{D}_1^+$]
The set $\hat{D}_1^+$ is not empty.
\end{proposition}
\begin{proof}
The set $\hat{D}_1^+$ is defined as the rectangle in $\bR^2$ bounded by the lines $\{x=k\}$, $\{x=\hat{u}_1\}$, $\{y=\hat{v}_1\}$ and $\{y=\hat{w}_1\}$.
It follows from \eqref{eq:u1} that $k < \hat{u}_1$, which in turn implies $\hat{w}_1 < \hat{v}$ if one compares the formulas \eqref{eq:w1} and \eqref{eq:v} that define these two constants.
\end{proof}

Now we are able to prove the main result of our paper.

\begin{theorem}[the absorbing set]\label{thm:absorbing}
The set $\hat{D}_1^+$ defined by \eqref{eq:d1plus} is an absorbing set for the Chialvo map~\eqref{eq:main}.
\end{theorem}

\begin{proof}
Take any $(x,y) \in \bR^n$. We shall prove that some iterate of this point will fall into $\hat{D}_1^+$ and will stay there forever.

By Proposition~\ref{prop:A}, the first iterate of $(x,y)$, as well as all the next iterates, are not in Region~$\hat{A}$ defined by \eqref{eq:regionA}.

If $(x,y)$ or any of its iterates are in Region~$\hat{B}$ defined by \eqref{eq:regionB} then by Proposition~\ref{prop:B} the trajectory of $(x,y)$ will leave $\hat{B}$ in a finite number of iterates and will not come back to $\hat{B}$. Therefore, the next iterates will be located to the right of the line $\{x = k\}$ and below the line $\{y = \hat{v}\}$, which is the union of Regions $\hat{C}_1$ and~$\hat{D}_1$ defined by \eqref{eq:regionC1} and \eqref{eq:regionD1}, respectively.

If $(x,y)$ or any of its iterates are in $\hat{C}_1 \cup \hat{D}_1$ then by Proposition~\ref{prop:C1} the next iterate is in Region~$\hat{D}_1$. By applying this proposition to the next iterates of $(x,y)$, we see that the entire forward trajectory from that point is contained in $\hat{D}_1$.

If $(x,y)$ or any of its iterates are in Region $\hat{E}_1 \subset \hat{D}_1$ then Proposition~\ref{prop:E1} combined with the just shown positive invariance of $\hat{D}_1$ imply that after a finite number of iterations the trajectory will reach $\hat{D}_1^+$ and will stay there forever. This completes the proof.
\end{proof}


\section{New construction of a weakly absorbing set for the Chialvo map}
\label{sec:weak}

Although the set $\hat{D}_1^+$ constructed in Section~\ref{sec:abs} is absorbing and thus solves the question posed in this paper, it might be tempting to construct a weakly absorbing set that might be smaller in some cases. The exponentiation in the formula for the Chialvo map~\eqref{eq:main} is challenging from the numerical point of view, so a tighter bound on $\hat{D}^+$ is of undeniable importance. Therefore, instead of taking the maximum of the right-hand side of the first equation of \eqref{eq:main} to obtain $\hat{u}_1$ in \eqref{eq:u1} that makes trajectories leave Region~$\hat{C}_1$ in just one iterate, in what follows we use the idea from \cite{courbage2010} and propose to determine a value of $\hat{u}_2$ that guarantees that trajectories ``move leftward'' in the corresponding Region~$\hat{C}_2$ and leave it after a finite number of iterates, like we did with Regions $\hat{B}$ and $\hat{E}_1$ (see Lemmas \ref{lem:downB} and~\ref{lem:upE1}, respectively). We also do not require that the new set $\hat{D}_2^+$ is positively invariant, but we still want to prove that every trajectory hits it in forward time in order to make sure it is weakly absorbing.


\subsection{\textbf{Region $\hat{C}_2$.}}
\label{sec:C2}
Assume that $y \leq \hat{v}$ and let us determine $\hat{u}_2 > k$ such that whenever $x \geq \hat{u}_2$ (note the weak inequality)  and $y \leq \hat{v}$, we have $\bar{x} < x$. For this purpose we need to determine a number $\hat{u}_2$ such that we can prove that
\begin{equation}
    \label{eq:C1}
    \bar{x} - x = x^2 \exp(y-x) + k - x < 0,
\end{equation}
for all $x \geq \hat{u}_2$, under the assumption that $y \leq \hat{v}$.

Inequality \eqref{eq:C1} can be written in the equivalent form as
\begin{equation}
    \label{eq:C1a}
    y < \ln \frac{x-k}{x^2} + x.
\end{equation}
However, we know that $y \leq \hat{v}$, so if we replace $y$ with $\hat{v}$  in Inequality~\eqref{eq:C1a} then Inequality~\eqref{eq:C1} will follow. We thus need to find $x_0$ large enough so that for every $x \geq x_0$ we have
\begin{equation}
    \label{eq:C2}
    \frac{c-bk}{1-a} < x + \ln \frac{x-k}{x^2}.
\end{equation}

\begin{lemma}
\label{lem:x0}
There exists a number $x_0$ such that Inequality~\eqref{eq:C2} is satisfied for every $x \geq x_0$. Moreover, such a number can be found effectively by means of a numerical algorithm.
\end{lemma}

\begin{proof}
Since $(x-k)/x^2 \approx 1/x$ for large $x$ and $\lim_{x \to \infty} x - \ln x = \infty$, there exists a number $x_0$ such that Inequality~\eqref{eq:C2} is satisfied for all $x \geq x_0$.
However, we are interested in finding such a number explicitly by means of an algorithm.
Moreover, we want to have this number as small as possible in order to obtain tight bounds for the weakly absorbing set $\hat{D}^+_2$ constructed in Section~\ref{sec:D2plus}.
In order to develop a reliable numerical procedure for finding a satisfactory number $x_0$, let us calculate the derivative of the right-hand side of Inequality \eqref{eq:C2}:
\begin{multline}
\label{eq:C3}
\frac{d}{dx}\left(x + \ln \frac{x-k}{x^2}\right) =
1 + \frac{x^2}{x-k} \, \frac{d}{dx}\left(\frac{1}{x} - \frac{k}{x^2}\right) \\
= 1 + \frac{x^2}{x-k}\left(-\frac{1}{x^2} + \frac{2k}{x^3}\right) =
\frac{x^2 - (k+1)x + 2k}{(x-k)x}.
\end{multline}
Notice that the denominator is positive if $x > k$. Solve the quadratic equation in the enumerator. If $(k+1)^2 - 8k \geq 0$, which happens for $k < 3 - 2 \sqrt{2} \approx 0.17$ and for $k > 3 + 2 \sqrt{2}$, then for all $x > x_1 := (k + 1 + \sqrt{k^2 - 6k + 1}) / 2$, the polynomial in the numerator is positive. In case $x_1 < k$, we must replace it with any number greater than $k$; we suggest $1.001 k$. Otherwise, if $(k+1)^2 - 8k < 0$, the numerator is always positive and we can take any number $x_1 > 0$ greater than $k$ (recall that $k > 0$); in order to keep some ``safety margin,'' we suggest to take $k$ increased by $0.1\%$ of $k$. We then begin a numerical procedure with $x = x_1$ and increase $x$ until we reach the value $x = x_0$ for which \eqref{eq:C2} is satisfied.
\end{proof}
Let us take
\begin{equation}
    \label{eq:u2}
    \hat{u}_2 := x_0 \text{ given by Lemma~\ref{lem:x0},}
\end{equation}
and define
\begin{equation}
\label{eq:regionC2}
\hat{C}_2 := \{\hat{u}_2 < x, \; y \leq \hat{v}\}.
\end{equation}
This construction implies that the following holds true:

\begin{lemma}[moving leftward in the closure of Region~$\hat{C}_2$]
\label{lem:leftC2}
Let $(x,y) \in \{\hat{u}_2 \leq x, \; y \leq \hat{v}\} \supset \hat{C}_2$. Then $\bar{x} < x$.
\end{lemma}

\begin{remark}
For the parameters in counterexample \eqref{eq:counter3}, we would have $\hat{u}_2 = 3.504$, which is indeed greater than the $x$ coordinate of the fixed point $(x,y)=(2.553,1.612)$ found for those parameter values in numerical simulations. Therefore, the fixed point is not in Region~$\hat{C}_2$.
\end{remark}

In contrary to what we could prove for Regions $\hat{B}$ (Lemma~\ref{lem:downB}) and~$\hat{E}_1$ (Lemma~\ref{lem:upE1}), Lemma~\ref{lem:leftC2} does not provide a lower bound on $x - \bar{x}$. Nevertheless, we have the following result.

\begin{proposition}[Region $\hat{C}_2$]
\label{prop:C2}
Every trajectory that starts in the set $\hat{C}_2$ leaves it in a finite number of iterations.
\end{proposition}

\begin{proof}
Suppose, for the sake of contradiction, that there is a trajectory $(x_i,y_i)_{i=0,1,\ldots}$ that stays in $\hat{C}_2$ forever. By Lemma~\ref{lem:leftC2}, the sequence $(x_i)_{i=0,1,\ldots}$ of the $x$ coordinates of the points in the trajectory is strongly monotonically decreasing ($x_{i+1} < x_{i}$ for all $i \geq 0$). We know by Theorem~\ref{thm:absorbing} that there exists $N > 0$ such that $(x_i,y_i) \in \hat{D}_1^+$ for all $i > N$. Note that $\hat{D}_1^+$ is a compact set. Then there exists a subsequence $(x_{k_j},y_{k_j})_{j=0,1,\ldots}$ of the trajectory $(x_i,y_i)_{i=0,1,\ldots}$ that is convergent to a point; denote it by $(x,y)$. Since $(x,y)$ is in the closure of $\hat{C}_2$, we can apply Lemma~\ref{lem:leftC2} to see that $\bar{x} < x$. Obviously, the subsequence $(x_{k_j+1},y_{k_j+1})_{j=0,1,\ldots}$ converges to $(\bar{x},\bar{y})$ with $\bar{x} < x$.
Consider $\tilde{x} := (x + \bar{x}) / 2$
and note that $x < \tilde{x} < \bar{x}$.
By convergence of the two subsequences, there exists a number $N_0$ such that
$x_{k_j} > \tilde{x}$ and $x_{k_j+1} < \tilde{x}$ for all $j > N_0$.
Take any $j_1 > N_0$ and $j_2 > j_1 + 1$.
Note that obviously $k_{j_1}+1 < k_{j_2}$.
Then $x_{k_{j_1}+1} < \tilde{x} < x_{k_{j_2}}$.
This contradicts the fact that $(x_i)_{i=0,1,\ldots}$ is decreasing, which completes the proof.
\end{proof}

Note that in the case of Region~$\hat{C}_2$, we do not have a result analogous to Lemma~\ref{lem:noB} for Region~$\hat{B}$ or Lemma~\ref{lem:noE1} for Region~$\hat{E}_1$. Proposition~\ref{prop:C2} does not rule out the situation that trajectories that left $\hat{C}_2$ come back to $\hat{C}_2$ after some iterations, and such trajectories can be indeed observed for some parameters; see the example discussed after Theorem~\ref{thm:weak}.


\subsection{\textbf{Regions $\hat{D}_2$ and $\hat{E}_2$.}}
\label{sec:E2}

Similarly to what we did in Sections \ref{sec:C1} and~\ref{sec:E1}, we use the number $\hat{u}_2$ instead of $\hat{u}_1$ to define the number $\hat{w}_1$ and Regions $\hat{D}_2$ and $\hat{E}_2$ in the following way:
\begin{eqnarray}
\label{eq:w2}
\hat{w}_2 & := & \frac{c - b \hat{u}_2}{1-a}, \\
\label{eq:regionD2}
\hat{D}_2 & := & \{k \leq x \leq \hat{u}_2,\; y \leq \hat{v}\}, \\
\label{eq:regionE2}
\hat{E}_2 & := & \{ k \leq x \leq \hat{u}_2,\; y < \hat{w}_2 \} \subset \hat{D}_2.
\end{eqnarray}
However, unlike it was in the case of Region $\hat{D}_1$, the set $\hat{D}_2$ is not positively invariant.
 Similarly, even though trajectories move upward in~$\hat{E}_2$, they do not need to stay in $\hat{D}_2$, and they can come back to~$\hat{C}_2$. Nevertheless, the following counterparts of Lemmas \ref{lem:upE1} and~\ref{lem:noE1} hold true (the proofs are the same as proofs of Lemmas \ref{lem:upE1} and~\ref{lem:noE1}, respectively):

\begin{lemma}[moving upward in Region $\hat{E}_2$]
\label{lem:upE2}
There exists a positive number $s > 0$ such that for every $(x,y) \in \hat{E}_2$ at least one of the following three conditions holds: $\bar{y} - y \geq s$, $(\bar{x},\bar{y}) \notin \hat{E}_2$, or $\bar{\bar{y}} - \bar{y} \geq s$.
\end{lemma}

\begin{lemma}[no return to Region $\hat{E}_2$]
\label{lem:noE2}
If $(x,y) \in \hat{D}_2 \setminus \hat{E}_2$ then $(\bar{x},\bar{y}) \notin \hat{E}_2$.
\end{lemma}


\subsection{\textbf{Region $\hat{D}_2^+$.}}
\label{sec:D2plus}
Similarly to \eqref{eq:d1plus}, we define Region $\hat{D}_2^+$ in the following way:
\begin{equation}
\label{eq:d2plus}
\hat{D}_2^+ := [k, \hat{u}_2] \times [\hat{w}_2,\hat{v}].
\end{equation}

In order to prove that Region~$\hat{D}_2^+$ is weakly absorbing, we will need an additional assumption.

\begin{theorem}[the weakly absorbing set]
\label{thm:weak}
Let us assume that
\begin{equation}
\label{eq:weak}
a+b \geq 1.
\end{equation}
Then the set $\hat{D}_2^+$ defined by \eqref{eq:d2plus} is a weakly absorbing set for the Chialvo map~\eqref{eq:main}.
\end{theorem}

Before we prove this theorem, let us justify the introduction of the additional assumption. For example, for the following values of the parameters:
\begin{equation}
a=0.5,\; b=0.01,\; c=0.98,\; k=0.9
\end{equation}
one can numerically find the following attracting orbit of period~$2$:
\begin{eqnarray}
(\bar{x},\bar{y}) & \approx & (4.01925,1.8954493) \in \hat{C}_2, \\
(x,y) & \approx & (2.8316796,1.8875321) \in \hat{E}_2.
\end{eqnarray}
The set $\hat{D}_2^+$ is bounded by the following values:
\begin{equation}
k=0.9,\; \hat{u}_2 \approx 3.4901791,\; \hat{w}_2 \approx 1.8901964,\; \hat{v} \approx 1.942.
\end{equation}
For the parameters in consideration, we have $\frac{b}{1-a} = 0.02 < 1$, and indeed the orbit oscillates between $\hat{E}_2$ and $\hat{C}_2$, and remains outside Region~$\hat{D}_2^+$.

Note, however, that the periodic orbit found is contained in $\hat{D}_1^+$, because $\hat{u}_1 \approx 4.6745998$ and $\hat{w}_1 \approx 1.866508$, so this is not a counterexample for Theorem~\ref{thm:absorbing}.

Before we prove Theorem~\ref{thm:weak}, we shall need the following result.

\begin{lemma}[no direct return from $\hat{E}_2$ to $\hat{C}_2$]
\label{lem:noEC}
Assume that $a+b \geq 1$. If $(x,y) \in \hat{E}_2$ then $(\bar{x},\bar{y}) := f(x,y) \in \hat{D}_2$.
\end{lemma}

\begin{proof}
Take any $(x,y) \in \hat{E}_2$.
By Proposition~\ref{prop:A}, $(\bar{x},\bar{y}) \notin \hat{A}$.
By Lemma~\ref{lem:noB}, $(\bar{x},\bar{y}) \notin \hat{B}$.
It remains to prove that $(\bar{x},\bar{y}) \notin \hat{C}_2$,
that is, $\bar{x} \leq \hat{u}_2$.

Let us first use the formula \eqref{eq:main} for $\bar{x}$, and then the inequality $y \leq \hat{w}_2 = \frac{c-b\hat{u}_2}{1-a}$ that defines the set $\hat{E}_2$ in~\eqref{eq:regionE2}.
\begin{equation}
\bar{x} =
x^2 \exp(y-x) + k \leq
x^2 \exp\big(\frac{c-b\hat{u}_2}{1-a} - x\big) + k \ldots
\end{equation}
Now we replace $-b\hat{u}_2$ with $-bk-b(\hat{u}_2-k)$ and we split the fraction into two:
\begin{equation}
\ldots = x^2 \exp \big(\frac{c - bk}{1-a} - \frac{b(\hat{u}_2-k)}{1-a} - x\big) + k \ldots
\end{equation}
Next, recall the inequality \eqref{eq:C2}. It holds true for all $x \geq \hat{u}_2$, in particular for $x = \hat{u}_2$. Let us use it to bound the first fraction under the exponent. As a consequence, we can continue as follows:
\begin{multline}
\ldots \leq x^2 \exp \big(\hat{u}_2 + \ln \frac{\hat{u}_2-k}{\hat{u}_2^2} - \frac{b(\hat{u}_2-k)}{1-a} - x\big) + k \\
= x^2 \exp \big(\ln \frac{\hat{u}_2-k}{\hat{u}_2^2} + (\hat{u}_2-x) - \frac{b}{1-a}(\hat{u}_2-k)\big) + k \\
= \frac{x^2}{\hat{u}_2^2}(\hat{u}_2-k) \exp \big((\hat{u}_2-x) - \frac{b}{1-a}(\hat{u}_2-k)\big) + k \ldots
\end{multline}
Since $x \leq \hat{u}_2$, we obviously have $\frac{x^2}{\hat{u}_2^2} \leq 1$. Also, $\hat{u}_2 - x \leq \hat{u}_2 - k$, because $x \geq k$. Recall the assumption that $a+b \geq 1$, which can be written as $\frac{b}{1-a} \geq 1$, and continue as follows:
\begin{equation}
\ldots \leq (\hat{u}_2-k) \exp \big((\hat{u}_2-k) - 1(\hat{u}_2-k)\big) + k = (\hat{u}_2-k) + k = \hat{u}_2,
\end{equation}
which is precisely the upper bound we need.
\end{proof}

\begin{proof}[Proof of Theorem~\ref{thm:weak}]
Take any $(x,y) \in \bR^n$. We shall prove that some iterate of this point will fall into $\hat{D}_2^+$.

By the same arguments as at the beginning of the proof of Theorem~\ref{thm:absorbing}, starting from some iterate, all the forward iterates of the point $(x,y)$ will be located in the union of Regions $\hat{C}_2$ and~$\hat{D}_2$ defined by \eqref{eq:regionC2} and \eqref{eq:regionD2}, respectively.

If the point $(x,y)$ or any its iterate is in Region~$\hat{C}_2$ then by Proposition~\ref{prop:C2}, after a finite number of iterations it will leave $\hat{C}_2$ and thus it will fall into $\hat{D}_2$. Lemma~\ref{lem:noEC} guarantees that the next iterations will remain in $\hat{D}_2$, and Lemma~\ref{lem:upE2} ensures that the trajectory will leave $\hat{E}_2$ in a finite number of steps, and thus fall into~$\hat{D}_2^+$. This completes the proof.
\end{proof}


\section{Numerical simulations}
\label{sec:num}

We conducted extensive numerical simulations in order to compare the sets $\hat{D}_1^+$ and $\hat{D}_2^+$ and to verify how tight the additional condition \eqref{eq:weak} in Theorem~\ref{thm:weak} is. Unless otherwise stated, we did the computations described here for all the values of
\begin{equation}
\label{eq:num}
\begin{array}{rcl}
a, b, c & \in & \{0.01, 0.02, \ldots, 0.99\} \subset (0,1), \\
k & \in & \{0.01, 0.02, \ldots, 0.99\} \cup \{1, 1.1, \ldots, 20\} \subset (0, \infty),
\end{array}
\end{equation}
which makes about $2.8 \cdot 10^{8}$ (over a quarter of a billion) combinations of the parameters.


\subsection{Advantage of $\hat{D}_2^+$ over $\hat{D}_1^+$ and practical construction of an absorbing set from $\hat{D}_2^+$.}
\label{sec:D1plusHuge}

In order to check the advantage of taking the weakly absorbing set $\hat{D}_2^+$ given by \eqref{eq:d2plus} over the absorbing set $\hat{D}_1^+$ given by \eqref{eq:d1plus}, we notice that $\hat{u}_1$ may be very large if $\hat{v}$ is large, which happens when $a \approx 1$ and $c-bk$ is large.
For example, for the following combination of the parameters:
\begin{equation}
\label{eq:u1badParam}
a = 0.9, \quad b = 0.1, \quad c = 0.9, \quad k = 0.1
\end{equation}
one obtains
\begin{equation}
\label{eq:u1badValues}
\hat{v} = 8.9, \quad \hat{u}_1 \approx 3969, \quad \hat{w}_1 \approx -3960.
\end{equation}
The region $\hat{D}_1^+$ resulting from these estimations is useless for any practical purposes; for example, subdividing it into squares of size $0.1$ for the purpose of numerical investigation would yield over $1.5$ billion ($1.5 \cdot 10^{9}$) grid elements.

However, the bounds computed using the method proposed in Section~\ref{sec:C2} provide a considerably better starting point for numerical investigation of dynamics, and the difference is profound. Indeed,
\begin{equation}
\label{eq:u2badValues}
\hat{v} = 8.9, \quad \hat{u}_2 \approx 11.34, \quad \hat{w}_2 \approx -2.337,
\end{equation}
which can be covered by some $1.26$ thousand ($1.26 \cdot 10^{3}$) squares of size $0.1$.

Let us notice that the values of the parameters \eqref{eq:u1badParam} in this example satisfy the assumption~\eqref{eq:weak} of Theorem~\ref{thm:weak}; therefore, the corresponding region $\hat{D}_2^+$ is weakly absorbing.

\begin{figure}[htbp]
\includegraphics[width=1.0\textwidth]{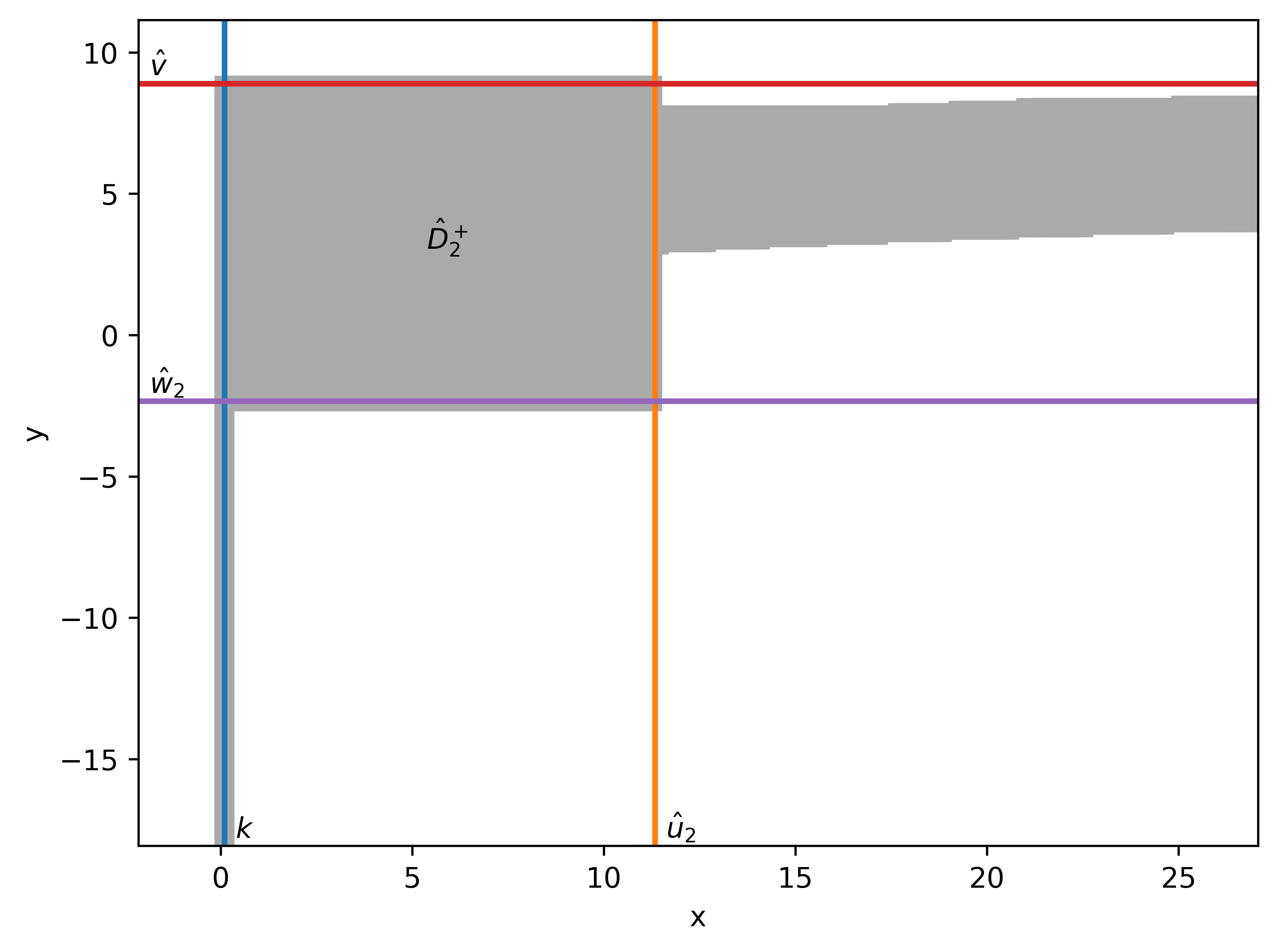}
\caption{\label{fig:tough}%
An outer approximation of the weakly absorbing set $\hat{D}_2^+$ and part of the absorbing set $\hat{P} := \bigcup_{n\geq 0} f^n(P)$, shown in grey, constructed for the Chialvo map with $a=0.9$, $b=0.01$, $c=0.9$ and $k=0.1$. The rectangles that form the set $\hat{P}$ are shown enlarged to make them clearly visible, also behind the lines. The ``arms'' of the set $\hat{P}$ actually extend downwards to $y \approx -430$ and rightwards almost to $x = 4400$.}
\end{figure}

In order to show a practical construction of an absorbing set from $\hat{D}_2^+$, we conducted the following experiment. We subdivided an outer numerical approximation $P := [0.099,11.35] \times [-2.34,8.91]$ of the set $\hat{D}_2^+$ into a uniform grid of $128 \times 128$ rectangles of the same size. Using interval arithmetic, we computed outer bounds for the images of these grid elements, and covered them by rectangles with respect to the same grid size. We added this outer approximation of $f(\hat{D}_2^+)$ to the set being constructed. We computed the images of the new grid elements, added them to the set, and we repeated this process until no further rectangles were added. In this way we constructed the smallest set $\hat{P}$, built of grid elements, containing $P$, and satisfying $f(\hat{P}) \subset \hat{P}$. This computation took about $1.5$ seconds on a modern laptop computer. The stabilization occurred already after two iterations, that is, $P \cup f(P) \cup f(f(P)) \subset \hat{P}$. The set $\hat{D}_2^+$ and part of the constructed set $\hat{P}$ are shown in Figure~\ref{fig:tough}.

The constructed set consisted of $681{,}774$ rectangles, and its area was about $5267$. Compare this to the area of $\hat{D}_2^+$, which was a little over $126$, and the area of $\hat{D}_1^+$, which was almost $1.6 \cdot 10^{7}$.
The constructed set was contained in the box $[0.099,4376.3] \times [-429.049,8.91]$. Note that part of the set was actually beyond $\hat{u}_1 \approx 3969$, but this should be attributed to overestimates coming from the rigorous numerical method applied. We remark that analogous computation in which the finer grid of $1024 \times 1024$ was used in $P$ resulted in the upper bound on the $x$ coordinate of some $4053$.


\subsection{Possibility of $\hat{D}_1^+ \subset \hat{D}_2^+$.}
\label{sec:D2subsetD1}

The features of $\hat{D}_2^+$ are weaker than those of $\hat{D}_1^+$, because the condition for $\hat{u}_2$ is less restrictive than that for $\hat{u}_1$. Therefore, one would expect that $\hat{D}_2^+ \subset \hat{D}_1^+$, like it was shown in the example in Section~\ref{sec:D1plusHuge}.
However, it turns out that in fact for some combinations of parameters one obtains $\hat{u}_1 < \hat{u}_2$. This implies that the absorbing set $\hat{D}_1^+$ is then contained in the weakly absorbing set $\hat{D}_2^+$. In such a situation there is no point to consider $\hat{D}_2^+$, because the bounds for the absorbing set $\hat{D}_1^+$ are tighter, while its features are stronger. We checked numerically the difference $\hat{u}_2 - \hat{u}_1$ among the parameters listed in~\eqref{eq:num}, and we found out that the largest difference encountered was only about $0.9$. It occurred for the following combination of parameters:
\begin{equation}
\label{eq:worstDiffParam}
a = 0.99, \; b = 0.99, \; c = 0.01, \; k = 0.04,
\end{equation}
which yielded
\begin{equation}
\label{eq:worstDiffValue}
\hat{u}_1 \approx 0.068, \; \hat{u}_2 \approx 0.956.
\end{equation}
In comparison to the advantage of using $\hat{u}_2$ instead of $\hat{u}_1$ shown in equations \eqref{eq:u1badValues} and~\eqref{eq:u2badValues}, we consider this difference negligible.


\subsection{Necessity of the assumption \eqref{eq:weak} in Theorem~\ref{thm:weak}}
\label{sec:evidence}

In order to verify whether the assumption \eqref{eq:weak} in Theorem~\ref{thm:weak} could be weakened, we checked for how many parameters out of those listed in \eqref{eq:num} that do not satisfy the inequality \eqref{eq:weak} we could really find a trajectory that does not reach Region~$\hat{D}_2^+$. Specifically, we conducted the following numerical simulation.

For every combination of the parameters listed in~\eqref{eq:num} that does not satisfy \eqref{eq:weak}, we first determined a rectangular region $R$ from which we would take starting points of trajectories to investigate. Since the set $\hat{D}_1^+$ is absorbing by Theorem~\ref{thm:absorbing}, a first attempt would be to take $R := \hat{D}_1^+$, as defined by~\ref{eq:d1plus}. However, as shown in Section~\ref{sec:D1plusHuge}, in some cases this region is huge, so we decided to take the region $4$ times larger than $\hat{D}_2^+$ in each direction if this happens, as shown in~\eqref{eq:R}. Since the maximum of $\hat{u}_2$ (that we observed in another numerical simulation) for all the parameters listed in \eqref{eq:num} was about $103.6$, we decided to take the threshold at $\hat{u}_1 = 105$. Specifically:
\begin{equation}
\label{eq:R}
R := \begin{cases}
[k,\; \hat{u}_1] \times [\hat{w}_1,\; \hat{v}] & \text{if } \hat{u}_1 < 105, \\
[k,\; \hat{u}_2+3(\hat{u}_2-k)] \times [\hat{w}_2-3(\hat{v}-\hat{w}_2),\; \hat{v}] & \text{otherwise.}
\end{cases}
\end{equation}

Next, we took a set of $10 \times 10$ initial conditions uniformly spread in~$R$. Then we checked if any forward iterate of each of these points was contained in $\hat{D}_2^+$. If it happened that within $1{,}000{,}000$ iterations the trajectory never crossed $\hat{D}_2^+$, it was reported as an argument in favor of the necessity of the assumption \eqref{eq:weak} in Theorem~\ref{thm:weak}. We emphasize the fact that we computed forward iterates even for initial conditions in $\hat{D}_2^+$ in case the images of these points were far away from~$\hat{D}_2^+$ and thus possibly missed when starting from the region~\eqref{eq:R}.

As a result, for $628{,}840$ out of $139{,}272{,}210$ parameter combinations ($0.45\%$) with $a+b < 1$, the program found a trajectory that did not enter $\hat{D}_2^+$, either because the trajectory stabilized with period-$2$ oscillations ($72.8\%$ of the cases), or the trajectory was still outside $\hat{D}_2^+$ after $1{,}000{,}000$ iterations ($27.2\%$ of the cases); note, however, that this was the case of the first encountered initial condition out of the $100 \times 100$ to be tested, so it is possible that another initial condition would provide the opposite result. In all these situations, the parameters were in the following ranges:
\begin{equation}
\label{eq:badRanges}
a \in [0.49, 0.98], \;
b \in [0.01, 0.07], \;
c \in [0.15, 0.99], \;
k \in [0.01, 20].
\end{equation}
It is an interesting observation that the set $\hat{D}_2^+$ was proved not to be weakly absorbing only for a half of the range of the parameter $a$, and only for very small values of the parameter $b$. This suggests that there might be some room for improvement in Theorem~\ref{thm:weak}.


\subsection{Number of iterates to reach $\hat{D}_1^+$}

Although we proved in Theorem~\ref{thm:absorbing} that every trajectory enters $\hat{D}_1^+$ in a finite number of steps, it turns out that this number can sometimes be very large. Indeed, for some combinations of parameters, the bounds on the ``speed'' along trajectories towards $\hat{D}_1^+$, denoted by $s$ in Lemmas \ref{lem:downB} and~\ref{lem:upE1}, can be very small. For example, for
\begin{equation}
a=0.99, \; b=0.77, \; c=0.01, \; k=0.48
\end{equation}
we have
\begin{equation}
\hat{v} \approx -35.96, \; \hat{u}_1 \approx 0.480000000048, \; \hat{w}_1 \approx -35.960000004,
\end{equation}
which yield
\begin{equation}
s(\hat{B}) \approx 9.75 \cdot 10^{-18}, \;
s(\hat{E}_1) \approx 1.44 \cdot 10^{-30}.
\end{equation}
Although the actual ``speed'' is much higher than these estimates, we still observed that for some points starting near $\hat{D}_1^+$ it takes over $100$ iterations to enter~$\hat{D}_1^+$.

Another problem is due to the exponential function in the first equation of \eqref{eq:main}. If the $y$ coordinate of some point is large and $x$ is small then the image of this point has a huge $x$ coordinate, which in turn yields extremely low $y$ coordinate of the next point. This situation is followed by an ascent of the trajectory according to $\bar{y} \approx a y$ which may be slow for $a \approx 1$. This scenario may lead to thousands of iterations until the trajectory reaches $\hat{D}_1^+$.
For example, for
\begin{equation}
a=0.99,\; b=0.03,\; c=0.06,\; k=0.41
\end{equation}
we have
\begin{equation}
\hat{v} \approx 4.77, \; \hat{u}_1 \approx 64.24, \; \hat{w}_1 \approx -187,
\end{equation}
and the point $(x_0,y_0) := (0.41, 64) \in \hat{D}_1^+$ is mapped to $(x_1,y_1) \approx (6.96 \cdot 10^{26}, 63.4)$, which is further mapped to $(x_2,y_2) \approx (0.41, -2.09 \cdot 10^{25})$, for which it takes over $5{,}200$ iterations to come back to $\hat{D}_1^+$.


\subsection{Tightness of the sets $\hat{D}_1^+$ and $\hat{D}_2^+$ in estimating the minimal absorbing set}

In order to check how tight our bounds in the form of the sets $\hat{D}_1^+$ and $\hat{D}_2^+$ are for estimating the actual global attractor of the system, we used the method introduced in \cite{Arai,chialvo2d} to compute a rigorous outer bound on the invariant part of $\hat{D}_1^+$ for a sample combination of parameter values considered in~\cite{Chialvo,chialvo2d}:
\begin{equation}
a=0.89,\; b=0.18,\; c=0.28,\; k=0.025.
\end{equation}
For these parameters, we have the following:
\begin{equation}
\hat{v} \approx 2.505,\;
\hat{u}_1 \approx 6.65, \;
\hat{u}_2 \approx 3.862, \;
\hat{w}_1 \approx -8.336, \;
\hat{w}_2 \approx -3.775.
\end{equation}
We plotted the set that was obtained when starting the computations from the rectangular region
\begin{equation}
P := [-0.637,7.313] \times [-9.42,3.589] \supset \hat{D}_1^+ \supset \hat{D}_2^+
\end{equation}

Specifically, the region $P$ was subdivided into a uniform grid of $256 \times 256$ rectangles. A directed graph that represented a rigorously computed outer approximation of the Chialvo map \eqref{eq:main} was created. Its strongly connected path components together with vertices on all the paths that joined these components were extracted. The set of rectangles corresponding to this subset of vertices was further subdivided into $2 \times 2$ rectangles each, and this computation was repeated. As justified in \cite{Arai}, the union of rectangles that correspond to the vertices of the final graph obtained in this way constitutes an outer bound for all the recurrent dynamics contained in $P$. The gradual refinements method applied here is also described and illustrated in~\cite{Arai}. We emphasize the fact that the computed estimation is rigorous, because it was obtained using interval arithmetic to calculate the outer approximation of the map on the grid elements.

\begin{figure}[htbp]
\includegraphics[width=1.0\textwidth]{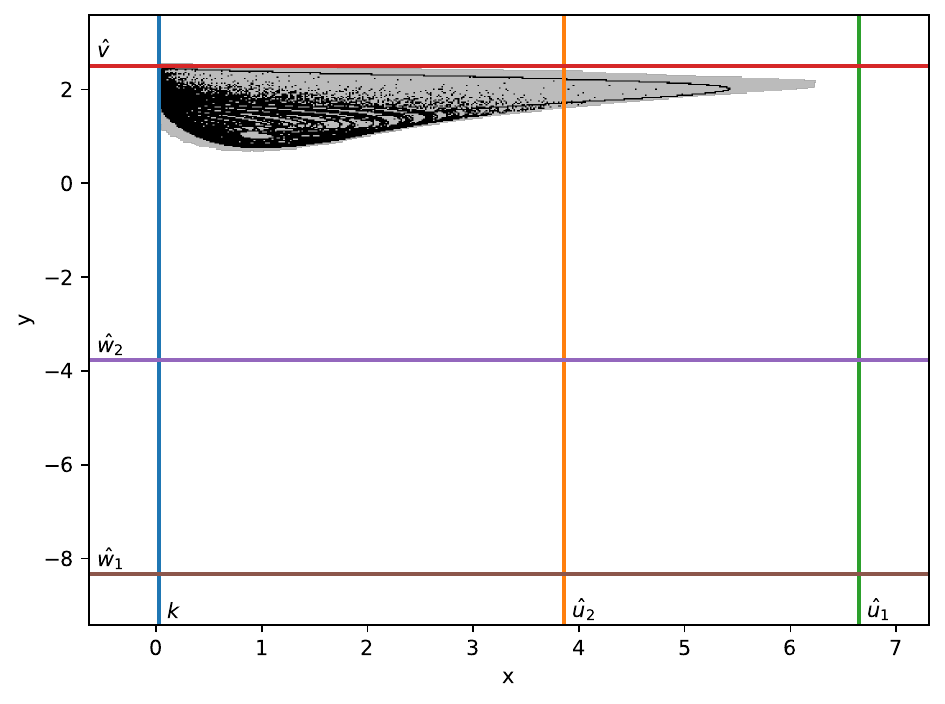}
\caption{\label{fig:sets}%
A rigorously computed outer bound for the absorbing set (gray), a numerical simulation of a trajectory on an attractor (black), and the constants that bound the regions $\hat{D}_1^+$ and~$\hat{D}_2^+$ for the Chialvo map with $a=0.89$, $b=0.18$, $c=0.28$ and $k=0.025$.}
\end{figure}

We additionally plotted a numerical approximation of the global attractor in the system. In order to obtain it, we started a trajectory in the center of $\hat{D}_1^+$, and after the initial $100{,}000$ iterations, we plotted the next $1{,}000{,}000$ points. The results are shown in Figure~\ref{fig:sets}.

The illustration shown in Figure~\ref{fig:sets} may serve as a summary of our discussion. The absorbing region $\hat{D}_1^+$ contains all the bounded recurrent dynamics in the system. The weakly absorbing region $\hat{D}_2^+$, although not necessarily absorbing in the sense of Definition~\ref{def:AbsorbingSet}, may serve as a starting point to reconstruct this dynamics. In fact, taking any of the two regions for conducting numerical analysis of the dynamics guarantees that no long-term dynamics will be missed. However, note that in our example the area of $\hat{D}_2^+$ is three times smaller than the area of $\hat{D}_1^+$, which provides certain advantage. As we showed in Section~\ref{sec:D1plusHuge}, for certain combinations of the parameters this difference may be profound. Therefore, considering the set $\hat{D}_2^+$ in addition to $\hat{D}_1^+$ is desired even though the properties of the former are weaker than those of the latter.


\subsection*{Acknowledgments}
We would like to express our gratitude to our colleague Justyna Signerska-Rynkowska for drawing our attention to the estimate $D^+$ for an absorbing set published in~\cite{courbage2010} that became a motivation for our work.

This research was supported by the National Science Centre, Poland, within the following grants:
Sheng~1 2018/30/Q/ST1/00228 (for G.G.),
and OPUS 2021/41/B/ST1/00405 (for P.P.).
Some of the computations were carried out using the computers of Centre of Informatics Tricity Academic Supercomputer \& Network.


\end{document}